\documentclass[12pt,a4paper]{article} 
\usepackage[T1]{fontenc}
\usepackage{libertine,libertinust1math} 
\usepackage{amsmath}
\usepackage{amsfonts}
\usepackage{amssymb}
\usepackage{ntheorem}
\usepackage{graphicx}
\usepackage{amsmath,amssymb}        
\usepackage{psfrag}
\usepackage{url}

\newtheorem{nth}{Theorem}

\def\RR{\hbox{I\kern-.2em\hbox{R}}}

\def\ds{\displaystyle}

\newcommand{\eqnsection}{
   \renewcommand{\theequation}{{\thesection.\arabic{equation}}}
   \makeatletter
   \csname @addtoreset\endcsname{equation}{section}
   \makeatother}
\eqnsection
\title{Existence and Uniqueness of BVPs Defined on Semi-Infinite Intervals:
Insight from the Iterative Transformation Method}
\author{Riccardo Fazio \\
Department of Mathematics, Computer Science,\\ Physical Sciences and Earth Sciences,\\
University of Messina \\
Viale F. Stagno D'Alcontres, 31 \\
98166 Messina, Italy \\
E-mail: \url{rfazio@unime.it} \\
Home-page: \url{http://mat521.unime.it/fazio}}
\date{\today}
\begin{document}               
\maketitle

\begin{abstract}
This work is concerned with the existence and uniqueness of boundary value problems defined on semi-infinite intervals.
These kinds of problems seldom admit exactly known solutions and, therefore, the theoretical information on their well-posedness is essential before attempting to derive an approximate solution by analytical or numerical means.
Our utmost contribution in this context is the definition of a numerical test for investigating the existence and uniqueness of solutions of boundary problems defined on semi-infinite intervals.
The main result is given by a theorem relating the existence and uniqueness question to the number of real zeros of a function implicitly defined within the formulation of the iterative transformation method. 
As a consequence, we can investigate the existence and uniqueness of solutions by studying the behaviour of that function.
Within such a context the numerical test is illustrated by two examples where we find meaningful numerical results.
\end{abstract}

\bigskip
\bigskip


\bigskip
\bigskip

\section{Introduction}
We consider the numerical solution of the boundary value problems (BVPs) defined on semi-infinite intervals belonging to the class of problems 
\begin{align}\label{eq:BVP:Class}
&{\ds \frac{d^{3}u}{dx^{3}}
= f\left(x, u, \frac{du}{dx}, \frac{d^2 u}{dx^2}\right)} \qquad ; \qquad   x\in [0, \infty) \nonumber \\[-1.5ex]
\\[-1.5ex]
&u(0) = u_0 \ , \quad {\ds \frac{du}{dx}(0) = v_0 \ , \qquad \frac{du}{dx}(x \rightarrow \infty) = v_\infty}  \ , \nonumber
\end{align}
where $ f(\cdot ,\cdot ,\cdot ,\cdot)$ satisfies appropriate smoothness conditions, $u_0$, $v_0$ and $v_\infty$ are given constants. 
BVPs 
belonging to (\ref{eq:BVP:Class}) arise mainly in boundary layer theory, see Schlichting and Gersten \cite{Schlichting:2000:BLT}, but also in fields of current interest like biology, chemistry, engineering, mechanics or physics. 
 
Before considering the application of any numerical method a question naturally arises: Is the BVP 
(\ref{eq:BVP:Class}) well-posed?
Indeed, few results are known about the existence and uniqueness question.
On this subject we can quote a paper by Weyl \cite{Weyl:1942:DES} about the celebrated Blasius problem \cite{Blasius:1908:GFK} and several papers concerning the more general Falkner-Skan model \cite{Falkner:1931:SAS}. 
Hence, any new idea for investigating that question is worth to be considered.

In this paper, we prove a theorem relating the existence and uniqueness question to the number of real zeros of a \lq \lq transformation function'' 
defined within the formulation of the iterative transformation method (ITM). 
Then, as a numerical test for the existence and uniqueness of solutions of (\ref{eq:BVP:Class}) we propose the application of the ITM 
to look for the number of real zeros of that function. 
In this context, the first application of that test was considered in \cite{Fazio:1991:ITM}.

The ITM discussed herein allows us to integrate numerically the class of BVPs (\ref{eq:BVP:Class}).
The method in point allows us to transform a BVP to a sequence of initial value problems. 
Such a transformation has also theoretical interest because it may represent an intermediate step to prove existence and uniqueness theorems (see \cite[pp. 7-13]{Keller}).

A first application of the ITM 
to a simple problem describing a biological model was given, on an intuitive basis, in \cite{Fazio:1990:SNA}. 
Moreover, other applications arose in connection with: the length determination for a tubular chemical reactor \cite{Fazio:1991:ITM}, the shock front propagation in rate-type materials \cite{Fazio:1992:MBH}, the classical Stefan's problem \cite{Fazio:1992:SAM}, and the spreading of a viscous fluid under the influence of gravity \cite{Fazio:1992:SAM}.
A further application of the ITM to the numerical solution of BVPs on infinite intervals governed by a third-order ordinary differential equation was described in \cite{Fazio:1996:NAN}.
A recent review concerning the applications of the ITM can be found, by the interested reader, in Fazio \cite{Fazio:2019:ITM}.

It was proved in \cite{Fazio:1994:NTM} that the ITM 
is an extension of the non-ITM 
proposed in \cite{Fazio:1990:SNA}.
The ITM 
can also be applied to two-point BVPs and in this case, it may be seen as an extension of the non-ITM 
discussed in \cite{Fazio:1990:NTM}. 
However, the applicability of non-ITMs is based on the invariance properties of the governing differential equation and of the given boundary conditions. 
As a consequence, non-ITMs are considered as \textit{ad hoc} methods (see \cite{Fox:1969:LBL}, \cite[pp. 35-36]{Meyer:1973:IVM}, \cite[p. 137]{Na:1979:CME}, \cite[p. 218]{Sachdev:1991:NOD}).
For the interested reader, a recent review reporting on the applications of the non-ITM can be found in Fazio \cite{Fazio:2019:NIT}.

To the best of our knowledge the ITM 
is the first transformation method applicable to a general class of problems.

The remainder of the paper is organized as follows. 
In section 2 we define the ITM. 
In section 3 we provide a simple constructive proof of the main result that relates the existence and uniqueness of solutions of BVPs defined on infinite intervals to the number of real zeros of the transformation function.
There, we point out how the present approach reduces the question of existence and uniqueness to find the number of real zeros of the transformation function.
A similar theorem has been proved by this author for a class of free BVPs governed by second-order differential equations, see \cite{Fazio:1997:NTE}.
Some guidelines on the definition of the transformation function and a discussion on its sensitivity analysis are developed in section 4. 
In section 5 we apply the iterative method to two BVPs 
describing, respectively, the Sakiadis problem and the Falkner-Skan model.
Our main concern is to verify numerically whether the considered problem has a unique solution.
In this context, the test indicates the existence and uniqueness of the solution for the first model and nonuniqueness for the second one. 
The last section is devoted to the final remarks and conclusions.
On the basis of the theorem proved in section 3, the proposed numerical test represents a possible way to investigate the existence and uniqueness of solution of BVPs.
Further evidence for this conclusion is given by the numerical experiments proposed in section~5 as well as by the study reported in \cite{Fazio:1991:ITM}.

In the following sections, we apply the ITM to the Sakiadis problem \cite{Sakiadis:1961:BLBa,Sakiadis:1961:BLBb}
\begin{align}\label{eq:S}
&{\displaystyle \frac{d^{3}f}{d\eta^3}} + \frac{1}{2} f {\displaystyle \frac{d^{2}f}{d\eta^2}} = 0 \nonumber \\[-1ex]
\\[-1ex]
&f(0) = 0 \ , \qquad {\displaystyle \frac{df}{d\eta}}(0) = 1 \ , \qquad
{\displaystyle \frac{df}{d\eta}}(\eta) \rightarrow 0 \quad \mbox{as}
\quad \eta \rightarrow \infty \ , \nonumber 
\end{align}

Moreover, we apply our ITM to the Falkner-Skan equation with relevant boundary conditions
\begin{align}\label{eq:FK}
& {\displaystyle \frac{d^{3}f}{d\eta^3}} + f 
{\displaystyle \frac{d^{2}f}{d\eta^2}} + \beta \left[ 1 - \left({\displaystyle
\frac{df}{d\eta}}\right)^2 \right] = 0 \nonumber \\[-1.2ex]
& \\[-1.2ex]
& f(0) = {\displaystyle \frac{df}{d\eta}}(0) = 0 \ , \qquad
{\displaystyle \frac{df}{d\eta}}(\eta) \rightarrow 1 \quad \mbox{as}
\quad \eta \rightarrow \infty \ , \nonumber
\end{align}
where $ f $ and $ \eta $ are appropriate similarity variables
and $ \beta $ is a parameter. 
This set is called the Falkner-Skan model, after the names of two English mathematicians
who first studied it \cite{Falkner:1931:SAS}.
Let us notice that for $\beta = 0$ the BVP (\ref{eq:FK}) reduces to the celebrated Blasius problem \cite{Blasius:1908:GFK}.
As pointed out by Na \cite[pp. 146-147]{Na:1979:CME}, if $ \beta \ne 0$, the BVP (\ref{eq:FK}) cannot be solved by a non-ITM.
Indeed, the governing differential equation in (\ref{eq:FK}) is not invariant with respect to any scaling group of point transformations.

The existence and uniqueness question for the problem (\ref{eq:FK}) is really a complex matter.
Assuming that $ \beta > 0 $ and under the restriction $ 0 < \frac{df}{d\eta} < 1 $, known as normal flow condition, Hartree \cite{Hartree:1937:EOF} and Stewartson \cite{Stewartson:1954:FSF} proved that the problem (\ref{eq:FK}) has a unique solution, whose first derivative tends to one exponentially.
Coppel \cite{Coppel:1960:DEB} and Craven and Peletier \cite{Craven:1972:USF} pointed out that the above restriction on the first derivative can be omitted when $ 0 \le \beta \le 1$.
For each value of the parameter $\beta$, there exists a physical solution with positive monotone decreasing second derivative, in the domain $[0, \infty)$, that approaches zero as the independent variables goes to infinity, as proved by Weyl \cite{Weyl:1942:DES}. 
In the case $ \beta > 1 $, the Falkner-Skan model loses the uniqueness and a hierarchy of solutions with reverse flow exists.
In fact, for $ \beta > 1$ Craven and Peletier \cite{Craven:1972:RFS} computed solutions for which $ \frac{df}{d\eta} < 0 $ for some value of $ \eta $.
In each of these solutions, the velocity approaches its limit exponentially in $\eta$. 
Here and in the following the term normal flow indicates that
the flow velocity has a unique direction, and instead, reverse flow means that the velocity
is both positive and negative in the integration interval. 

The considered problem has also multiple solutions for $\beta_{\min} < \beta < 0$, as reported by  Veldman and van de Vooren \cite{Veldman:1980:GFS}, with the minimum value of $\beta$ given by 
\begin{equation}\label{eq:bmin}
\beta_{\min} = -0.1988\dots \ .
\end{equation} 
In this range there exist two physical solutions, one for normal flow and one
for reverse flow.  
For $\beta = \beta_{\min}$ only one solution exists.
Finally, for $\beta < \beta_{\min}$ the problem has no solution at all. 
Our interest here is to apply the ITM to the range of $\beta$ where multiple solutions are admitted, in particular, we would like to compute numerically the famous solutions of Stewartson \cite{Stewartson:1954:FSF,Stewartson:1964:TLB}.
As we shall see in the following sections the obtained results are in very good agreement with those available in the literature.

Hartree \cite{Hartree:1937:EOF} was the first to solve numerically the Falkner-Skan model.
Then, Cebeci and Keller \cite{Cebeci:1971:SPS} applied shooting and parallel shooting methods requiring asymptotic boundary condition to be imposed at a changing unknown boundary in the computation process. 
As a result, they reported convergence difficulties, which can be avoided by moving towards more complex methods. 
Moreover, to guarantee reasonable accuracy, they were forced to use a {\it small enough} step-size and extensive computations for the solution of the IVPs.
Na \cite[pp. 280-286]{Na:1979:CME} described the application of invariant imbedding to the Falkner-Skan model.
For this problem, a modified shooting method \cite{Asaithambi:1997:NMS} and finite-difference methods \cite{Asaithambi:1998:FDM,Asaithambi:2004:SOF} were presented by Asaithambi.
Sher and Yakhot \cite{Sher:2001:NAS} defined a new approach to solve this problem by shooting from infinity, using some simple analysis of the asymptotic behaviour of the solution at infinity.
Kuo \cite{Kuo:2003:ADT} used a differential transformation method, to compute a series solution of the Falkner-Skan equation.
Asaithambi \cite{Asaithambi:2005:SFE} proposed a faster shooting method by using a recursive evaluation of Taylor coefficients.
Zhang and Chen \cite{Zhang:2009:IMS} investigated a modification of the shooting method, where the computation of the Jacobian matrix was obtained by solving two IVPs. 
A Galerkin-Laguerre spectral method was defined and applied to the Falkner-Skan model by Auteri and Quartapelle \cite{Auteri:2012:GLS}.

\section{The iterative transformation method}
By requiring the invariance of (\ref{eq:BVP:Class}) with respect to a given transformation group we characterize a subclass of problems (see \cite{Fazio:1990:SNA} for the stretching and the spiral group and \cite{Fazio:1990:NVT} for the translation group).
As a consequence non-ITMs are applicable only to special classes of problems. 
To overcome this drawback we let the function $f(\cdot , \cdot , \cdot , \cdot )$  in (\ref{eq:BVP:Class}) to depend also on a numerical parameter $ h $ and we use $h$ also to modify the initial conditions to get invariant initial conditions, and in this way, we introduce the problem
\begin{align}\label{eq:BVP:Class:h}
&{\ds \frac{d^{3}u}{dx^{3}}
= F\left(x, u, \frac{du}{dx}, \frac{d^2 u}{dx^2}, h \right)} \qquad ; \qquad   x\in [0, \infty) \nonumber \\[-1.5ex]
\\[-1.5ex]
&u(0) = h^{1/\sigma} u_0 \ , \quad {\ds \frac{du}{dx}(0) = h^{(1-\delta)/\sigma} v_0} \ , \qquad
{\ds \frac{du}{dx}}(x \rightarrow \infty) = v_\infty   \ . \nonumber
\end{align}
This allows us to consider the extended stretching group
\begin{equation}\label{eq:stretching}
x^{*} =  \lambda ^{\delta} x \quad ;  
\quad u^{*} =  \lambda  u \quad ;
\quad h^{*} =  \lambda ^{\sigma} h 
\end{equation}
where $ \delta $ and $ \sigma  $ are constants different from zero and $\lambda$ is the exponential of the group parameter, and therefore it is not zero. 
Our intention is to require the invariance of the governing differential equation in (\ref{eq:BVP:Class:h}) with respect to (\ref{eq:stretching}).
To this end, we have to require that the following relation 
\begin{equation}\label{eq:Fversusf}
F\left(x , u,\frac{du}{dx}, \frac{d^2 u}{dx^2}, 1\right) = h^{(1-3\delta)/\sigma} f\left(h^{-\delta/\sigma} x, h^{-1/\sigma} u, h^{(\delta-1)/\sigma)}\frac{du}{dx}, h^{(2\delta-1)/\sigma)}\frac{d^2 u}{dx^2}\right) \ ,
\end{equation}
holds true.
We impose (\ref{eq:Fversusf}) in order to recover (\ref{eq:BVP:Class}) from (\ref{eq:BVP:Class:h}) as $ h \rightarrow 1 $.
Therefore, (\ref{eq:BVP:Class}) can be embedded in the BVP
\begin{align}\label{eq:BVP:Class:Embedded}
&{\ds \frac{d^{3}u}{dx^{3}} = h^{(1-3\delta)/\sigma} f\left(h^{-\delta/\sigma} x, h^{-1/\sigma} u, h^{(\delta-1)/\sigma)}\frac{du}{dx}, h^{(2\delta-1)/\sigma)}\frac{d^2 u}{dx^2}\right)} \nonumber \\
\\[-1.5ex]
&u(0) = h^{1/\sigma} u_0 \ , \quad {\ds \frac{du}{dx}(0) = h^{(1-\delta)/\sigma} v_0} \ , \qquad
{\ds \frac{du}{dx}}(x \rightarrow \infty) = v_\infty  \ . \nonumber
\end{align}
Of course, a given problem belonging to (\ref{eq:BVP:Class}) can always be embedded into (\ref{eq:BVP:Class:Embedded}) for every values of $ \delta $ and $ \sigma $
different from zero.
The requested invariance means that (\ref{eq:stretching}) transforms a problem of form (\ref{eq:BVP:Class:Embedded}) to a problem of the same form but with different values of $h$ and $v(\infty)$ if $v(\infty) \ne 0$.

From a numerical viewpoint, we have to consider the initial value problem (IVP)
\ in the starred variables defined by
\begin{align}\label{eq:IVP*}
&{\ds \frac{d^{3}u^{*}}{dx^{*3}} = h^{*(1-3\delta )/\sigma }
f\left(h^{*-\delta /\sigma } x^{*}, h^{*-1/\sigma
} u^{*}, h^{*(\delta -1)/\sigma } \frac{du^{*}}{dx^{*}}, h^{*(2\delta -1)/\sigma }\frac{d^{}u^{*}}{dx^{*2}}\right)} \nonumber \\[-1.5ex]
\\[-1.5ex]
&u^*(0) = h^{1/\sigma} u_0 \ , \quad {\ds \frac{du^*}{dx^*}(0) = h^{(1-\delta)/\sigma} v_0} \ , \nonumber
\end{align}
where $ \delta $ and $ \sigma $ have to be considered as fixed values. 
If for every value of $ h^{*} $ the problem (\ref{eq:IVP*}) is well posed, then $ v^{*}(0) $ is uniquely defined. 
Consequently, the application of invariance considerations allows us to obtain, if $v_\infty \ne 0$,
\begin{equation}\label{eq:lambda}
\lambda = \left(\frac{\frac{du^*}{dx^*}(\infty)}{v_{\infty}}\right)^{1/(1-\delta)}.
\end{equation}  
If $v_\infty = 0 $, then we can use $h$ to get a non-invariant boundary condition at infinity.
In fact, we can use the new side condition
\begin{equation}\label{eq:newBC}
\frac{du}{dx}(\infty) = 1 - h^{(1-\delta)/\sigma} \ .
\end{equation}
Therefore, if $v_\infty = 0$ then we apply the new boundary condition (\ref{eq:newBC}) and get $\lambda$ by
\begin{equation}\label{eq:lambda2}
\lambda = \left({\ds \frac{du^*}{dx^*}(\infty) + {h^*}^{(1-\delta)/\sigma}}\right)^{1/(1-\delta)}.
\end{equation}  
In both cases, for a positive value of $ \lambda $ we get the value of $h$ given by
\begin{equation}\label{eq:h}
h = \lambda ^{-\sigma } h^{*} \ .
\end{equation}  
A solution of (\ref{eq:BVP:Class}) is determined when the value of $ h=1 $ is obtained from (\ref{eq:h}).
If $\delta $ and $ \sigma $ are fixed, then $ \lambda  $ is a function of $ h^{*} $ only, i.e., $ \lambda  = \lambda (h^{*}) $ where, in general, the functional form of $ \lambda (\cdot ) $ is not known.
Thus, we are interested in the zeros of the \lq \lq transformation function''
\begin{equation}\label{eq:Gamma}
\Gamma (h^{*}) = [\lambda (h^{*})]^{-\sigma } h^{*} - 1 \qquad ; \qquad \lambda > 0 \ . 
\end{equation}  
This function is defined implicitly by the solution of the IVP (\ref{eq:IVP*}).

Along the lines of the analysis sketched above the ITM can be defined as follows:
the original BVP is embedded into (\ref{eq:BVP:Class:h}) by fixing non-vanishing values of $ \delta $ and $ \sigma $;
 starting with suitable values of $ h^{*}_{0} $ and $ h^{*}_{1} $ a root-finder method, secant, regula-falsi, quasi-Newton or bisection, is used to define a sequence $ h^{*}_{j}, $ for $ i = 2, 3, \dots, \  $.
At each iteration, a related sequence given by $ \Gamma (h^{*}_{j}), $ for $ i = 0, 1, 2, \dots, \  $ is defined according to (\ref{eq:Gamma});
 suitable termination criteria have to be used to verify whether $ \Gamma (h^{*}_{j}) \rightarrow  0 $ as $ j \rightarrow \infty$;
 a numerical approximation of $ u(x) $ can be obtained by rescaling the numerical solution $u^*(x^*)$ corresponding to $h = 1$.

\section{Existence and uniqueness}
The main result of this paper can be stated as follows: for a given BVP 
defined on a semi-infinite interval the existence and uniqueness question is reduced to find the number of real zeros of the transformation function. 
This result is proved below.

\begin{nth}\label{Th:exuni}
Let us assume that $ \delta $ and $ \sigma $ are fixed and 
that for every value of $ h^{*} $ the IVP (\ref{eq:IVP*}) is well-posed. 
Then, the BVP (\ref{eq:BVP:Class}) has a unique solution if and only if the transformation function has a unique real zero; nonexistence (nonuniqueness) of the solution of (\ref{eq:BVP:Class}) is equivalent to nonexistence of real zeros (existence of more than one real zero) of $ \Gamma (\cdot ) $.
\end{nth}
{\bf Proof by invariance considerations.} We begin by proving that there exists a one-to-one and onto correspondence between the set of solutions of (1.1) and the set of real zeros of the transformation function.
The thesis is an evident consequence of this result.

The mentioned correspondence can be defined as follows.
First, for every values of $ \delta $ and $ \sigma $ different from zero, given a solution $ u(x)$ of (\ref{eq:BVP:Class:h}) we can associate to it the real zero of $ \Gamma (\cdot ) $ defined by 
$ h^{*} = {v_\infty}/{\frac{du^*}{dx^*}(\infty)} $. 
The related value of $ \lambda = {h^*}^{1/\sigma} $ allows us to verify by substitution in (\ref{eq:Gamma}) that we have defined a real zero of $ \Gamma (\cdot ) $. 
Moreover, because of the invariance with respect to (\ref{eq:stretching}) the related solution of (\ref{eq:IVP*}) is given by 
$ u^{*}(x^{*}) = h^{*1/\sigma } u(h^{* \delta/\sigma} x) $.
 
Second, according to the definition of the transformation function, in general to each real zero $ h^{*} $ of $ \Gamma (\cdot ) $ there is related a solution $ u^{*}(x^{*}) $ of the IVP (\ref{eq:IVP*}). 
The condition for $u^{*}(x^{*}), h^{*} $ to be transformed by (\ref{eq:stretching}) to $u(x), 1 $ (where $ u(x) $ is defined on $ [0, \infty) $) is $ \lambda = h^{*1/\sigma} $. 
Due to $ \lambda = h^{*1/\sigma} $, we have that $ \frac{du^*}{dx^*}(\infty) = h^{*(1-\delta)/\sigma} v_\infty $, $ u^*(0) = h^{*1/\sigma} u_0 $ and $ \frac{du^*}{dx^*} (0) =h^{*(1-\delta)/\sigma} v_0 $, so that the relation (\ref{eq:lambda}) implies that $ u(x) $ verifies the boundary condition at zero and at infinity given in (\ref{eq:BVP:Class}). 
Hence, to each real zero of $ \Gamma (\cdot ) $ we can associate a solution of (\ref{eq:BVP:Class}).
Again $ \lambda = h^{*1/\sigma } $, so that $ u(x) = h^{*-1/\sigma } u^{*}(h^{*-\delta/\sigma} x^{*}) $.

It is easily seen that both a right and a left inverse of our correspondence are fixed by means of the relations defined above. 
Therefore, the correspondence is one-to-one and onto.
Finally, it follows that if one of the two sets is empty then the other one is empty too.$ \hfill \square $

\noindent
{\bf Remarks:} Before proceeding further some remarks are in order. 
First of all, as far as IVPs are concerned, the theory of well-posed problems is developed in detail in several classical books (see, for instance, Hartman \cite[Chapters 2, 3 and 5]{Hartman:1982:ODE}).
In particular, the continuous dependence of the solution on parameters holds true provided suitable regularity conditions on $ f(\cdot, \cdot, \cdot , \cdot ) $ are fulfilled.
Second, if for every value of $ h^{*} $ we assume $ \lambda (h^{*}) > 0 $, then for $ \delta \neq 0 $ and each fixed value of $ h^* $ the scaling 
$ [\lambda (h^{*})]^{-\delta } x^{*} : [0, \infty) \rightarrow [0, \infty) $ is one-to-one and onto whereas the function of $ h^* $ defined by
$ [\lambda (h^{*})]^{-\sigma } h^{*} : \RR \rightarrow \RR $ may not be one-to-one for $ \sigma \neq 0 $. 
Therefore, since $ \Gamma (h^{*}) = [\lambda (h^{*})]^{-\sigma } h^{*} - 1 $, the transformation function may not be one-to-one.
Third, by studying the behaviour of the transformation function it is possible to get useful insight into the existence and uniqueness question. 
 
\section{The transformation function}
In this section, we provide some guidelines for the definition of $ \Gamma (\cdot). $
Let us define the function
\[
 \Xi (h^{*}, \delta , \sigma) = 
[\lambda (h^{*}, \delta , 
\sigma)]^{-\sigma } h^{*} - 1 \  ,
\]
by considering $\delta $ and $ \sigma $ as variables.
Of course, for different values of these variables, we obtain different transformation functions.
In any case, the analysis proposed in the previous section does not depend upon the choice of the values of $\delta $ and $ \sigma $.
A sensitivity analysis for the dependence of $ \Xi (\cdot, \cdot, \cdot) $ with respect to $\delta $ and $ \sigma $ is not necessary because the ITM is defined for fixed values of these variables.
As a consequence, once we have fixed the values of $\delta $ and $ \sigma $ it is unimportant if we use only approximations of these values.

Moreover, by the proof of the theorem above we know in advance that every real zero of the transformation function belongs to $ \RR ^{+} $ and therefore we can restrict the domain of $ \Gamma (\cdot ) $ to $ \RR ^{+} $ whereupon its range will be reduced to $ (-1, \infty ) $.
We remark that the values of $ \delta $ and $ \sigma $ can always be chosen to consider only positive values of $ h^{*} $.

On the other hand, the sensitivity of $ \Gamma (\cdot ) $ (and consequently of $ \Xi (\cdot, \cdot, \cdot) $) with respect to $ h^* $ is of paramount interest.
In fact, the numerical determination of the roots of $ \Gamma (\cdot) = 0 $ is a well-conditioned problem if and only if 
$ | {d \Gamma}/{d h^*} | $ at the root is bigger than the machine rounding unit (a simple proof of this statement is possible under suitable regularity conditions on $ \Gamma (\cdot) $).
Therefore, it is relevant to verify the sensitivity of $ \Gamma (\cdot) $ at each zero.
Of course, since $ \Gamma (\cdot) $ is not given explicitly we cannot evaluate its sensitivity directly.
Hence, monitoring of the sensitivity will require an increase in the computational cost.
A simple procedure is to apply as root-finder the secant method because in this case, we compute at each iteration a finite difference approximation for the first derivative of $ \Gamma (\cdot) $.
Some caution has to be used because finite difference approximations of derivatives are prone to rounding errors that, for really small increments, could dominate the approximations.

Let us consider the definition of $ \lambda $ in (\ref{eq:lambda}) or in (\ref{eq:lambda2}), since the initial conditions are invariant with respect to the extended scaling group (\ref{eq:stretching}), $ \lambda (h^{*}, \delta , \sigma) $ must have the following functional form
\[
\lambda (h^{*}, \delta , \sigma)  = h^{*1/\sigma}
\lambda_{0} \bigl(\delta , \sigma \bigr)
\]
so that
\[
 \Xi (h^{*}, \delta , \sigma) =
[\lambda_{0} \bigl(\delta , 
\sigma \bigr)]^{-\sigma } - 1 \ \ \ .
\]
In other words, $ \Gamma(\cdot) $ depends on $ h^{*-\delta/\sigma} $, instead of $ h^{*} $ only.
This can be used when the zeros of the transformation function are very close to one endpoint of its domain.
In particular, we can modify the value of $ \delta $, $ \sigma $ or both values to define a $\Gamma(\cdot)$ function with more amenable zeros.

\section{Numerical tests and results}
In this section, we test the ITM 
with respect to the existence and uniqueness question and obtain a numerical solution for the given BVP.

\subsection{Sakiadis problem}
In order to apply the ITM to (\ref{eq:S}) we introduce, see Fazio \cite{Fazio:2015:ITM}, the extended problem
\begin{align}\label{eq:Sakiadis:mod}
&{\displaystyle \frac{d^3 f}{d \eta^3}} + \frac{1}{2} f
{\displaystyle \frac{d^{2}f}{d\eta^2}} = 0 \nonumber \\[-1.5ex]
& \\[-1.5ex]
&f(0) = 0 \ , \qquad {\displaystyle \frac{df}{d\eta}}(0) = h^{1/2} \ , \qquad
{\displaystyle \frac{df}{d\eta}}(\eta \rightarrow \infty) \rightarrow 1 - h^{1/2} \ . \nonumber
\end{align}
In (\ref{eq:Sakiadis:mod}), the governing differential equation and the two initial conditions are invariant, whereas the asymptotic boundary condition is not invariant, with respect to the extended scaling group
\begin{equation}\label{eq:scaling}
f^* = \lambda f \ , \qquad \eta^* = \lambda^{-1} \eta \ , \qquad 
h^* = \lambda^{4} h \ .   
\end{equation}
Moreover, it is worth noticing that the extended problem (\ref{eq:Sakiadis:mod}) reduces to the Sakiadis problem  (\ref{eq:S}) for $h=1$.
This means that in order to find a solution of the Sakiadis problem we have to find a zero of the transformation function 
\begin{equation}\label{eq:AA2.9}
\Gamma (h^{*}) = \lambda^{-4} h^* - 1 \ , 
\end{equation}  
where the group parameter $ \lambda $ is defined by the formula
\begin{equation}\label{eq:lambda3}
\lambda = \left[\displaystyle \frac{df^*}{d\eta^*}(\eta_\infty^*)+{h^*}^{1/2}\right]^{1/2} \ ,
\end{equation}
and to this end we can use a root-finder method.

Let us notice that $\lambda$ and the transformation function are defined implicitly by the solution of the IVP
\begin{align}\label{eq:Sakiadis:IVP}
&{\displaystyle \frac{d^3 f^*}{d \eta^{*3}}} + \frac{1}{2} f^*
{\displaystyle \frac{d^{2}f^*}{d\eta^{*2}}} = 0 \nonumber \\[-1.5ex]
& \\[-1.5ex]
&f^*(0) = 0 \ , \quad {\displaystyle \frac{df^*}{d\eta^*}}(0) = {h^*}^{1/2}, \quad
{\displaystyle \frac{d^2f^*}{d\eta^{*2}}}(0) = \pm 1  \nonumber \ .
\end{align}
In particular, we are interested to compute $\frac{df^*}{d\eta^*}(\eta_\infty^*)$, where $\eta_\infty^*$ is a suitable truncated boundary, an approximation of the asymptotic value $\frac{df^*}{d\eta^*}(\infty)$, which is used in the definition (\ref{eq:lambda2}) of $\lambda$. 

It is evident that our numerical method is based on the behaviour of the transformation function.
Our interest is to study the behaviour of this function with respect to its independent variable.
We notice that because of the two terms $h^{*1/2}$, which have been introduced in the modified boundary conditions in (\ref{eq:Sakiadis:mod}), we are allowed to consider only positive values of $h^*$. 

\begin{figure}[!hbt]
	\centering
\psfrag{h*}[][]{$ h^* $} 
\psfrag{G}[][]{$ \Gamma(h^*) $} 
\psfrag{df}[][]{${\displaystyle \frac{d^2f^*}{{d\eta^*}^2}} (0) = 1$} 
\includegraphics[width=.8\textwidth]{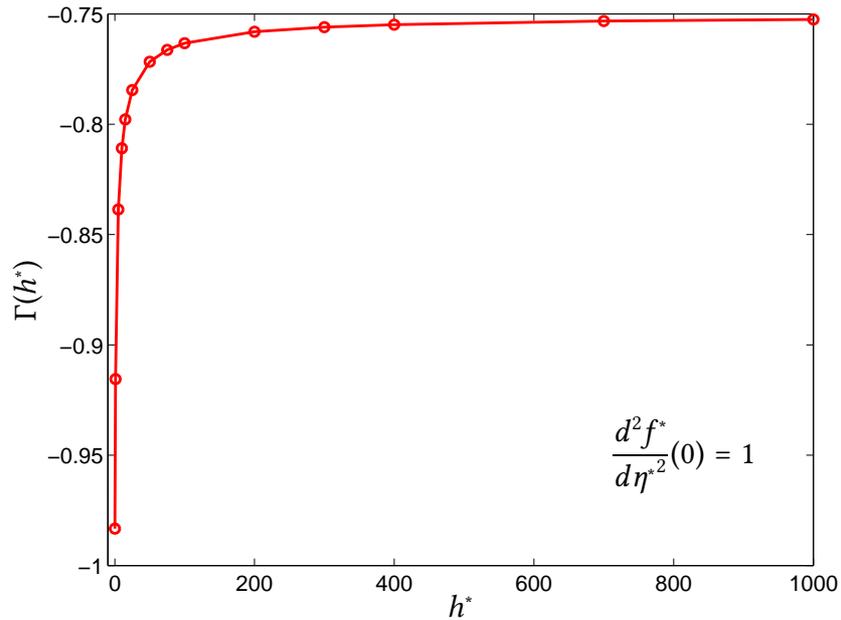} 
\caption{Plot of $\Gamma(h^*)$ for ${\displaystyle \frac{d^2f^*}{{d\eta^*}^2}}(0) = 1$.} 
	\label{fig:Gamma1}
\end{figure}

\begin{figure}[!hbt]
	\centering
\psfrag{h*}[][]{$ h^* $} 
\psfrag{G}[][]{$ \Gamma(h^*) $} 
\psfrag{df}[r][]{${\displaystyle \frac{d^2f^*}{{d\eta^*}^2}} (0) = -1$} 
\includegraphics[width=.9\textwidth]{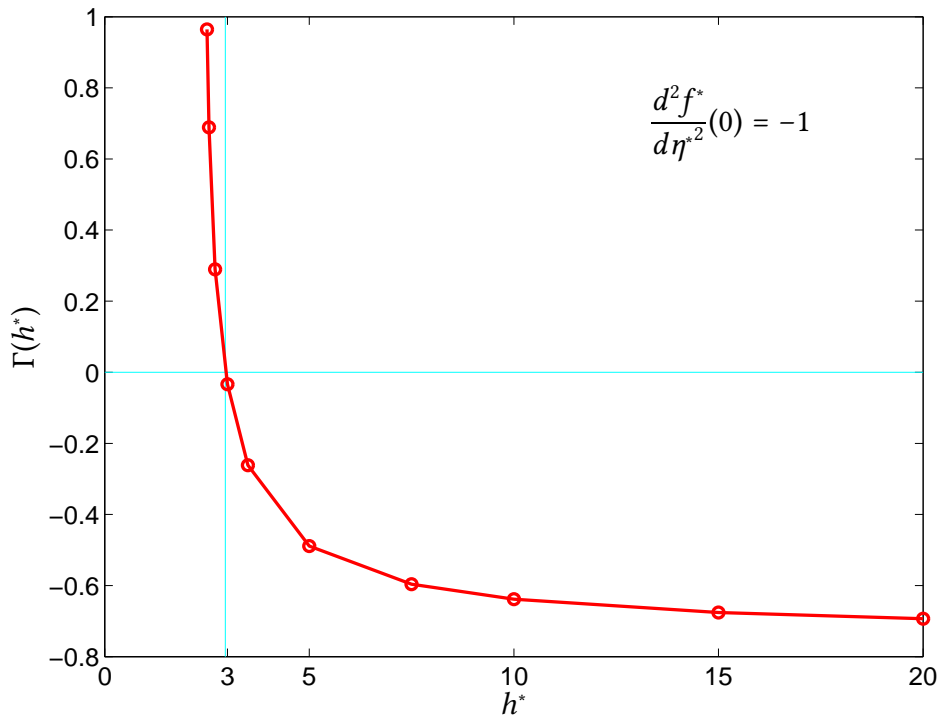} 
\caption{Plot of $\Gamma(h^*)$ for ${\displaystyle \frac{d^2f^*}{{d\eta^*}^2}} (0) = -1$.} 
	\label{fig:Gamma2}
\end{figure}

From our numerical study concerning the dependence of $\Gamma$ with respect to the missing initial condition $\frac{d^2f^*}{{d\eta^*}^2} (0)$ we have used the results plotted on figures \ref{fig:Gamma1}-\ref{fig:Gamma2}.
Each o-symbol represents a numerical solution of the IVP (\ref{eq:Sakiadis:IVP}) with the corresponding value of $h^*$. 
The solid line joining these symbols is used to show the behaviour of the transformation function.  
By considering the results reported in figure \ref{fig:Gamma1} we realize that the missing initial condition cannot be positive.

For a negative missing initial condition the numerical results are shown on figure \ref{fig:Gamma2}.
It is evident from figure \ref{fig:Gamma2} that the transformation function has only one zero and, by the theorem proved above, this means that the considered problem has one and only one solution.
Moreover, we remark that the tangent to the $\Gamma$ function at its unique zero and the $h^*$ axis define a large angle.
From a numerical viewpoint, this means that the quest for the $h^*$ corresponding to $h=1$ is a well-conditioned problem.

As far as the numerical results reported in this section are concerned, the ITM was applied by setting the truncated boundary $ \eta_\infty^* = 10$.
Moreover, these results were obtained by an adaptive fourth-order Runge-Kutta IVP solver. 
The adaptive solver uses a relative and absolute error tolerance, for each component of the numerical solution, both equal to $1 \mbox{D}{-06}$. 
Here and in the following the notation $\mbox{D}-k = 10^{-k}$ means double-precision arithmetic.

The initial value solver was coupled with the simple secant root-finder with a convergence criterion given by
\begin{equation}\label{eq:conv}
\left|\Gamma(h^*)\right| \le 1 \mbox{D}{-09} \ .
\end{equation}

\begin{figure}[!hbt]
	\centering
\psfrag{e}[][]{$ \eta $} 
\psfrag{f}[][]{$f$} 
\psfrag{df}[][]{$ {\displaystyle \frac{df}{d\eta}} $} 
\psfrag{ddf}[][]{$ {\displaystyle \frac{d^2f}{d\eta^2}} $} 
\includegraphics[width=.9\textwidth]{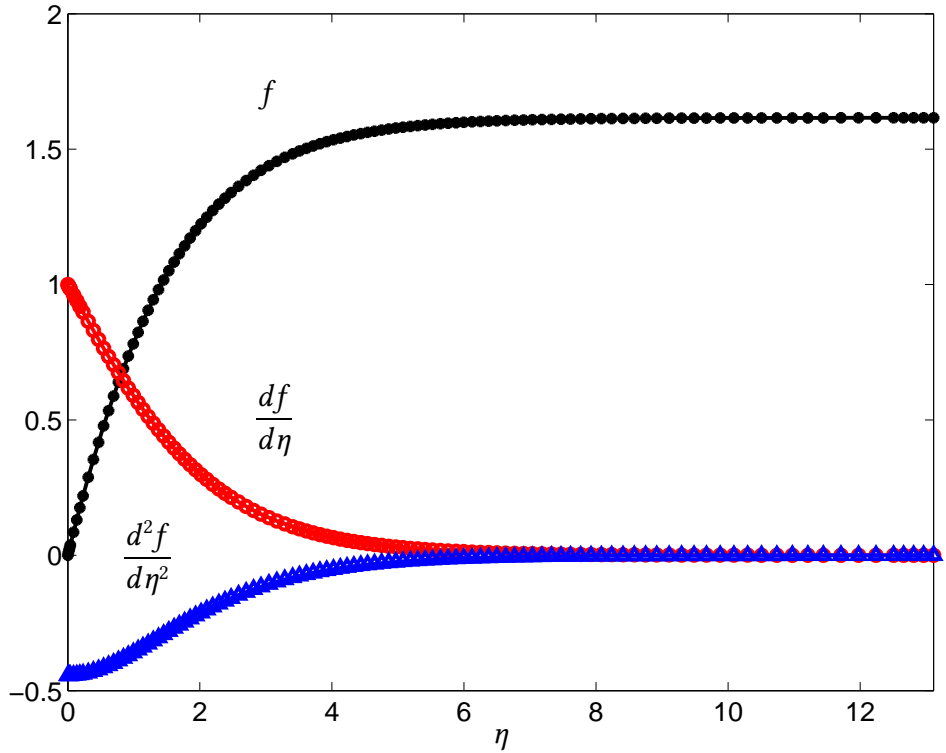} 
\caption{Sakiadis solution via the ITM.} 
	\label{fig:Sakiadis}
\end{figure}

The implementation of the secant method is straightforward.
The only difficulty we have to face is related to the choice of the initial iterates.
In this context, the study of the transformation function of figure \ref{fig:Gamma2} was really helpful.
In table \ref{tab:itera} we list the iterations of our ITM.
\begin{table}[!htb]
\caption{Iterations and numerical results: secant root-finder.}
\begin{center}{ \renewcommand\arraystretch{1.3}
\begin{tabular}{rr@{.}lcr@{.}lr@{.}l} 
\hline%
{$j$}
&
\multicolumn{2}{l}%
{$h_j^*$} 
& {$\lambda_j$}
& \multicolumn{2}{c}%
{$\Gamma(h_j^*)$} 
& \multicolumn{2}{c}%
{${\displaystyle \frac{d^2f}{{d\eta}^2}} (0)$} \\[1.5ex]
\hline%
$0$  & $2$ & $5$       & $1.061732$ & $0$  & $967343$ & $-0$ & $835517$ \\
$1$  & $3$ & $5$       & $1.475487$ & $-0$ & $261541$ & $-0$ & $311310$ \\
$2$ & $3$ & $287172$ & $1.417981$ & $-0$ & $186906$ & $-0$& $350743$\\
$3$ & $2$ & $754191$ & $1.229206$ & $0$ & $206411$ & $-0$ & $538426$ \\
$4$ & $3$ & $033897$ & $1.339089$ & $-0$ & $056455$ & $-0$ & $416458$ \\
$5$ & $2$ & $973826$ & $1.318081$ & $-0$ & $014749$ & $-0$ & $436690$ \\
$6$ & $2$ & $952581$ & $1.310382$ & $0$ & $001407$ & $-0$ & $444433$ \\
$7$ & $2$ & $954432$ & $1.311058$ & $-3$ & $23\mbox{D}{-05}$ & $-0$ & $443745$ \\
$8$ & $2$ & $954391$ & $1.311043$ & $-6$ & $93\mbox{D}{-08}$ & $-0$ & $443761$ \\
$9$ & $2$ & $954391$ & $1.311043$ & $3$ & $42\mbox{D}{-12}$ & $-0$ & $443761$ \\
\hline
\end{tabular}}
\label{tab:itera}
\end{center}
\end{table}
The last iteration of table \ref{tab:itera} defines our numerical approximation that is shown, for the reader convenience, on figure \ref{fig:Sakiadis}.
This solution was computed by rescaling, note that by rescaling we get ${\eta^*_\infty} < {\eta}_{\infty}$, and this means that we have reduced the computational cost, where the chosen truncated boundary was ${\eta^*_\infty} = 10$ in our case.

\subsection{The Falkner-Skan model}
In order to apply an ITM to (\ref{eq:FK}) we have to embed it to an extended model and require the invariance of this last model with respect to an extended scaling group of point transformations.
This can be done in several ways that are all equivalent.
In fact, the modified model can be written as 
\begin{align}\label{eq:abf2}
& {\displaystyle \frac{d^{3}f}{d\eta^3}} + f 
{\displaystyle \frac{d^{2}f}{d\eta^2}} + \beta \left[ h^{4/\sigma} - \left({\displaystyle
\frac{df}{d\eta}}\right)^2 \right] = 0 \ , \nonumber \\[-1.2ex]
& \\[-1.2ex]
& f(0) = {\displaystyle \frac{df}{d\eta}}(0) = 0 \ , \qquad
{\displaystyle \frac{df}{d\eta}}(\eta) \rightarrow 1 \quad \mbox{as}
\quad \eta \rightarrow \infty \ , \nonumber
\end{align}
and the related extended scaling group is given by
\begin{equation}\label{eq:scalinv}
f^* = \lambda f \ , \qquad \eta^* = \lambda^{-1} \eta \ , \qquad 
h^* = \lambda^{\sigma} h \ ,   
\end{equation}
where $\sigma$ is a parameter.
In the following we set $\sigma = 4$; for the choice $\sigma=8$ see \cite{Fazio:1994:FSE}.
In \cite{Fazio:1994:FSE}, a free boundary formulation of the Falkner-Skan model was considered and
numerical results were computed for the Homann flow ($ \beta = 1/2 $) as well as for the Hiemenz flow ($ \beta = 1 $).

From a numerical point of view the request to evaluate 
$ \frac{d f}{d \eta} (\infty) $ cannot be fulfilled.
Several strategies have been proposed in order to provide an approximation of this value.
The simplest and widely used one is to introduce, instead of infinity, a suitable truncated boundary.
For the sake of simplicity we apply, following T\"opfer \cite{Topfer:1912:BAB}, some preliminary computational tests to find a suitable value for the truncated boundary.

\begin{figure}[!htb]
	\centering
\psfrag{h*}[][]{$h^*$} 
\psfrag{G}[][]{$\Gamma(h^*)$} 
\includegraphics[width=.7\textwidth]{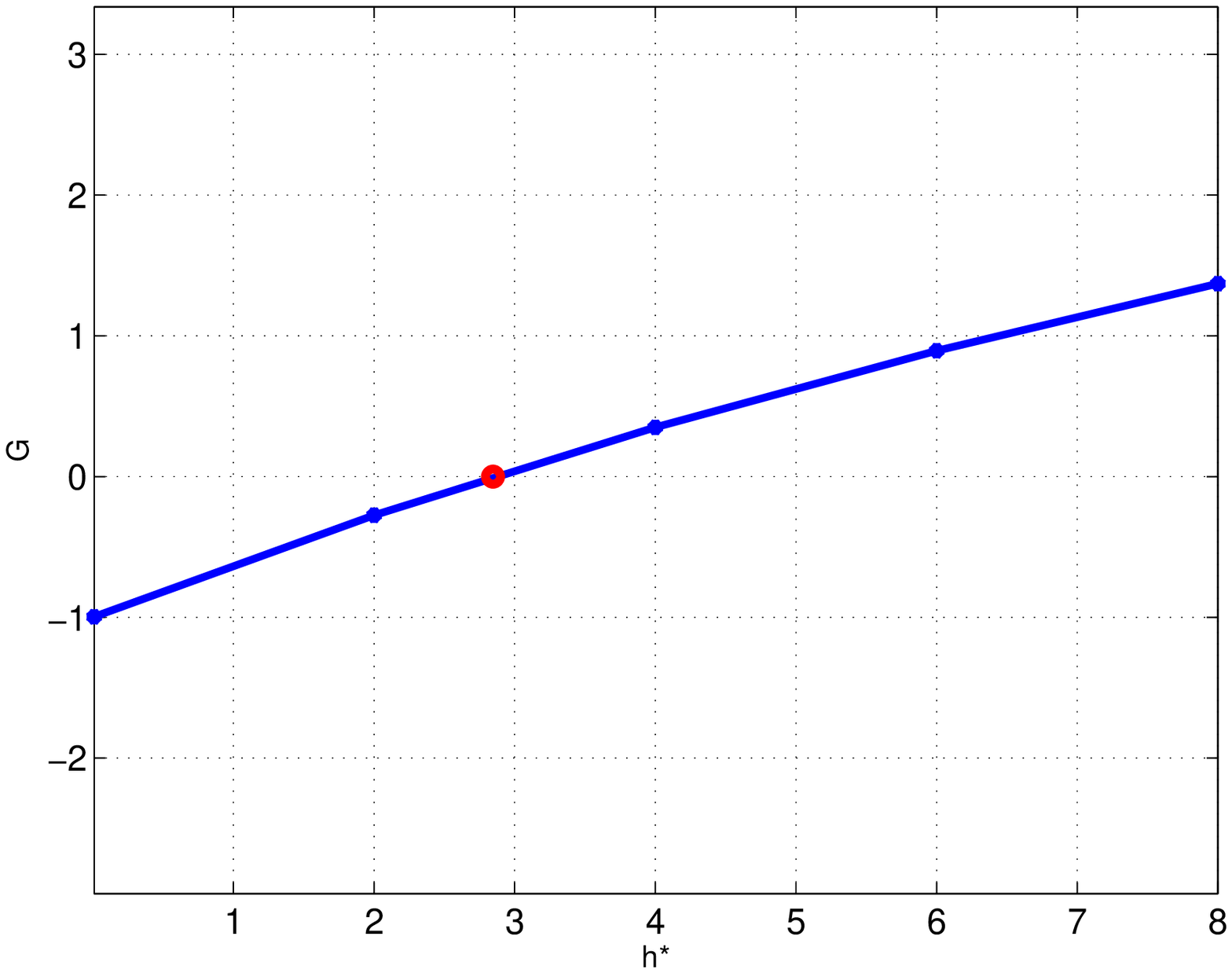} \\
\includegraphics[width=.7\textwidth]{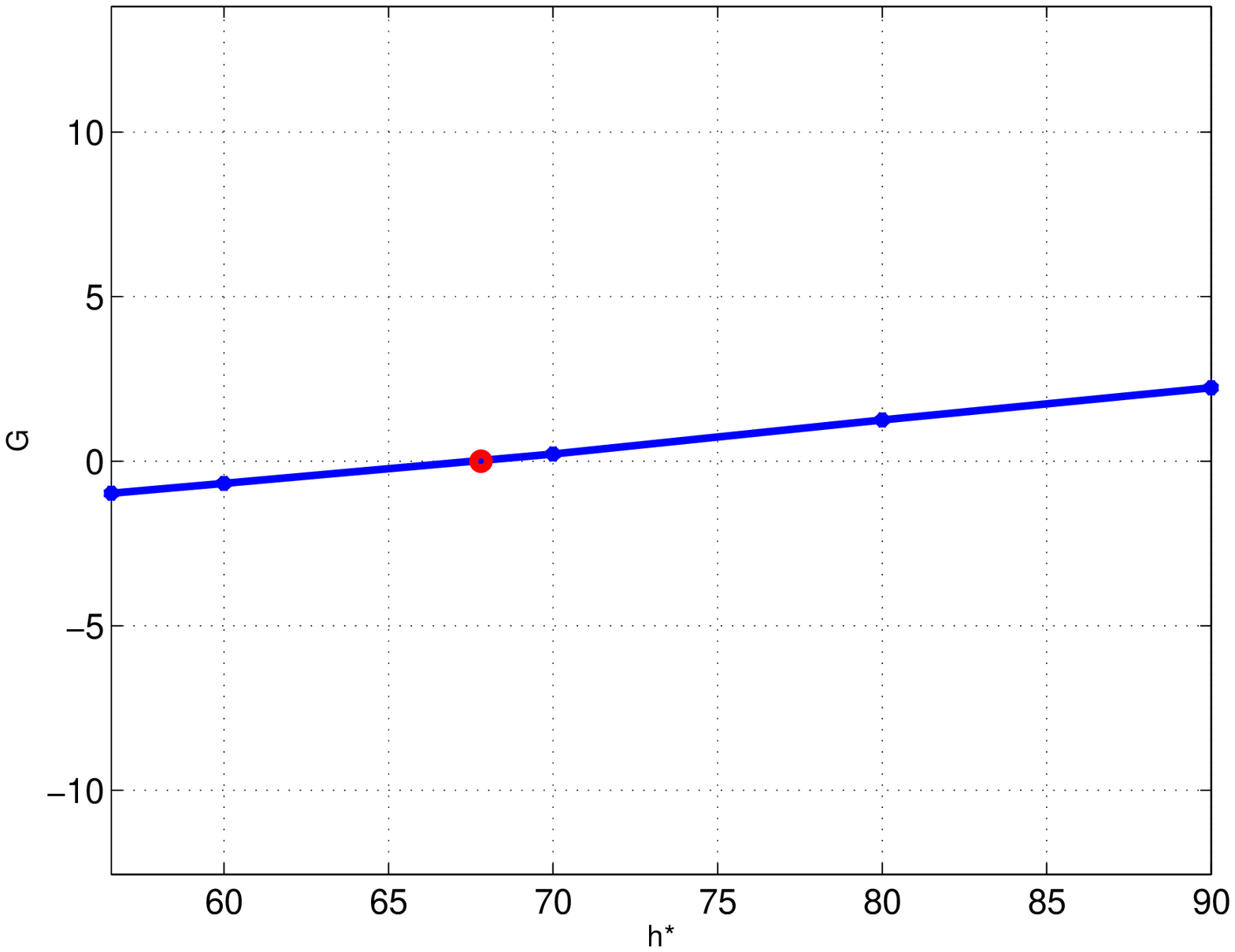}
\caption{Two cases of the $\Gamma(h^*)$ function: top and bottom frames are related to normal and reverse flow solutions, respectively.} 
	\label{fig:FSHGhstar}
\end{figure}

At each iteration of the ITM, we have to solve the IVP 
\begin{eqnarray}\label{eq:IVP2}
& {\displaystyle \frac{d^{3}f^*}{d\eta^{*3}}} + f^* 
{\displaystyle \frac{d^{2}f^*}{d\eta^{*2}}} + \beta \left[ h_j^{*} - \left({\displaystyle
\frac{df^*}{d\eta^*}}\right)^2 \right] = 0 \nonumber \\[-1.2ex]
& \\[-1.2ex]
& f^*(0) = {\displaystyle \frac{df^*}{d\eta^*}}(0) = 0 \ , \qquad
{\displaystyle \frac{d^2f^*}{d\eta^{*2}}}(0) = \pm 1  \ . \nonumber
\end{eqnarray}
Tables~\ref{tab:Itera1} and~\ref{tab:Itera2} list the numerical iterations obtained for a sample value of $\beta$.
We notice that we solve an IVP governed by a different differential equation for each iteration because the Falkner-Skan equation is not invariant under every scaling group of point transformations.
As a first case we have chosen $ \beta=-0.01$.
The data listed in tables~\ref{tab:Itera1} and~\ref{tab:Itera2} have been obtained by solving the extended Falkner-Skan model on $ \eta^* \in [0, 20] $ by setting 
\[
 \frac{d^2f^*}{d\eta^{*2}}(0) = \pm 1 \ ,
\]
respectively.
\begin{table}[!htb]
\caption{Iterations for $\beta=-0.01$ with ${\ds \frac{d^2f^*}{d\eta^{*2}}(0)} = 1$. Here and in the following the $\mbox{D}-k = 10^{-k}$ means a double precision arithmetic.}
\vspace{.5cm}
\renewcommand\arraystretch{1.3}
	\centering
		\begin{tabular}{cr@{.}lr@{.}lcr@{.}l}
\hline 
{$j$} &
\multicolumn{2}{c}%
{$h_j^*$}
& \multicolumn{2}{c}%
{$\Gamma(h_j^*)$}
&
{${\ds \frac{|h_j^*-h_{j-1}^*|}{|h_j^*|}}$}
& \multicolumn{2}{c}%
{${\displaystyle \frac{d^2f}{d\eta^2}(0)}$} \\[1.2ex]
\hline
0 &  5 &            &     0 & 631459 & & 0 & 431723 \\
1 & 10 &            &      1 & 791425  & & 0 & 384034 \\
2 &  2 & 278111 & $-$0 & 182888 & 3.389602 & 0 & 454658 \\
3 &  2 & 993420 &     0 & 0465208 & 0.238960 & 0 & 454658 \\
4 &  2 & 848366 &     9 & 5$\mbox{D}-04$ & 0.050925 & 0 & 456418 \\
5 &  2 & 845340 & $-$5 & 0$\mbox{D}-06$ & 0.001064 & 0 & 456455 \\
6 &  2 & 845356 &     6 & 1$\mbox{D}-08$ & 5.6$\mbox{D}-06$ & 0 & 456455 \\
7 &  2 & 845355 &     7 & 3$\mbox{D}-10$ & 6.7$\mbox{D}-08$ & 0 & 456455 \\
\hline			
		\end{tabular}
	\label{tab:Itera1}
\end{table}
\begin{table}[!htb]
\caption{Iterations for $\beta=-0.01$ with ${\ds \frac{d^2f^*}{d\eta^{*2}}(0)} = -1$.}
\vspace{.5cm}
\renewcommand\arraystretch{1.3}
	\centering
		\begin{tabular}{cr@{.}lr@{.}lcr@{.}l}
\hline 
{$j$} &
\multicolumn{2}{c}%
{$h_j^*$}
& \multicolumn{2}{c}%
{$\Gamma(h_j^*)$}
&
{${\ds \frac{|h_j^*-h_{j-1}^*|}{|h_j^*|}}$}
& \multicolumn{2}{c}%
{${\displaystyle \frac{d^2f}{d\eta^2}(0)}$} \\[1.2ex]
\hline
0 &  75 &            &     0 & 731890 & & $-$0 & 059237 \\
1 & 150 &           &     5 & 263092 & & $-$0 & 092368 \\
2 &  62 & 885833 & $-$0 & 443040 & 1.385275 & $-$0 & 028870 \\
3 &  69 & 649620 &     0 & 181067 & 0.097112 & $-$0 & 046991 \\
4 & 67 & 687299 & $-$0 & 011297  & 0.028991 & $-$0 & 042016 \\
5 & 67 & 802542 & $-$2 & 1$\mbox{D}-04$ & 0.001700 & $-$0 & 042315 \\
6 & 67 & 804749 & 2 & 8$\mbox{D}-07$      & 3.3$\mbox{D}-05$ & $-$0 & 042321 \\
7 & 67 & 804746 & 7 & 9$\mbox{D}-10$    & 4.3$\mbox{D}-08$ & $-$0 & 042321 \\
\hline			
		\end{tabular}
	\label{tab:Itera2}
\end{table}
In both cases, we achieved convergence of the numerical results within seven iterations.
Let us now investigate the behaviour of the transformation function.
Figure~\ref{fig:FSHGhstar} shows $\Gamma(h^*)$ with respect to $h^*$ for the two cases reported in these tables.
The unique zero of the transformation function is marked by a circle.
The same scale has been used for both axes.
As it is easily seen, in both cases, we have a monotone increasing function.
We notice on the top frame, corresponding to normal flow, that the tangent to the $\Gamma$ function at its unique zero and the $h^*$ axis define a large angle.
This is important from a numerical viewpoint because in such a case we face a well-conditioned problem.
On the other hand, this is not the case for the function plotted on the bottom frame of the same figure. 
The meaning is clear, reverse flow solutions are more challenging to compute than normal flow ones.

\begin{figure}[!htb]
	\centering
\psfrag{e}[][]{\small $ \eta $} 
\psfrag{f}[][]{\small ${\displaystyle f(\eta), \frac{df}{d\eta}, \frac{d^2f}{d\eta^2}}$} 
\psfrag{fe}[][]{$ {\displaystyle \frac{df}{d\eta}} $} 
\psfrag{fee}[][]{$ {\displaystyle \frac{d^2f}{d\eta^2}} $} 
\includegraphics[width=.7\textwidth]{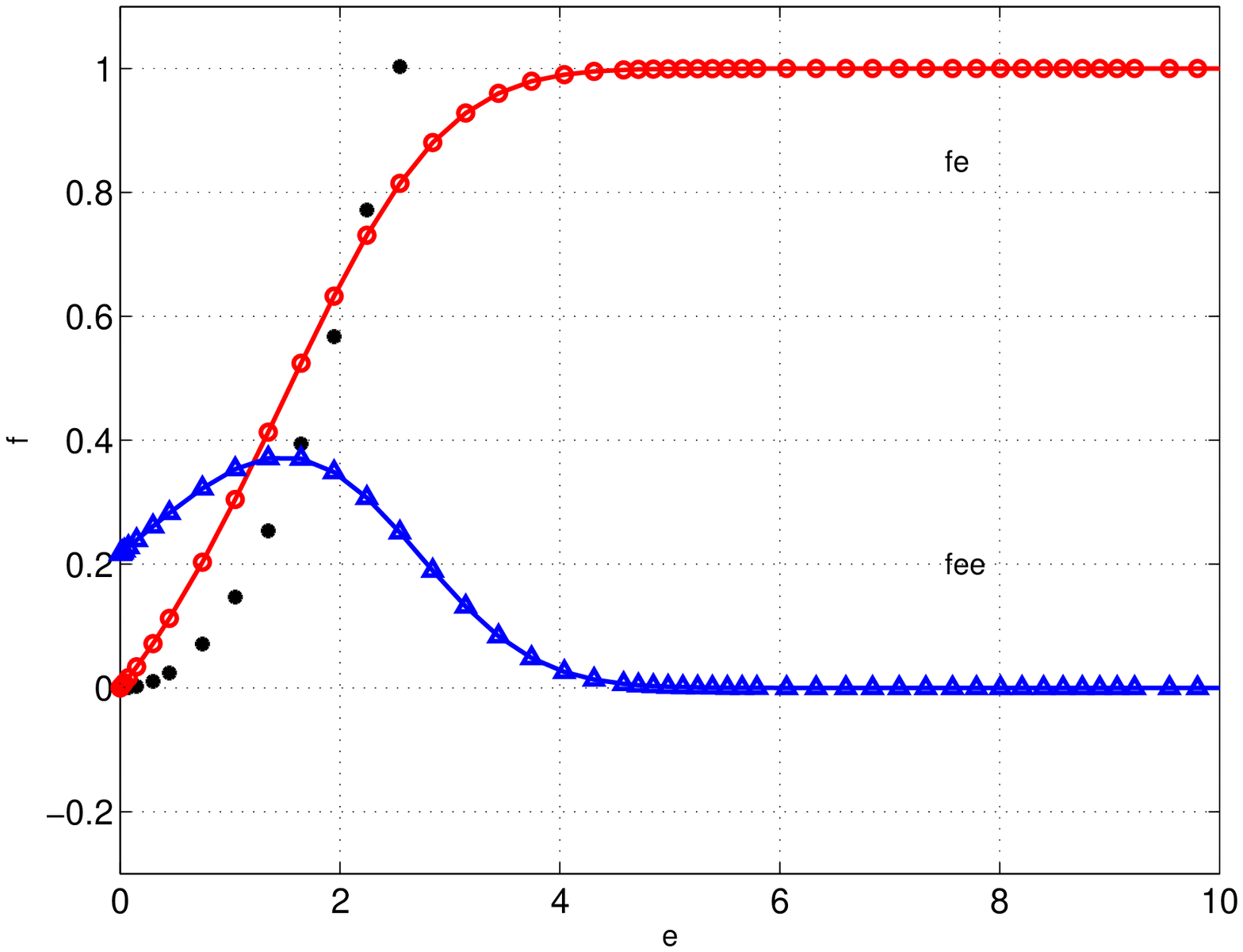} \\
\includegraphics[width=.7\textwidth]{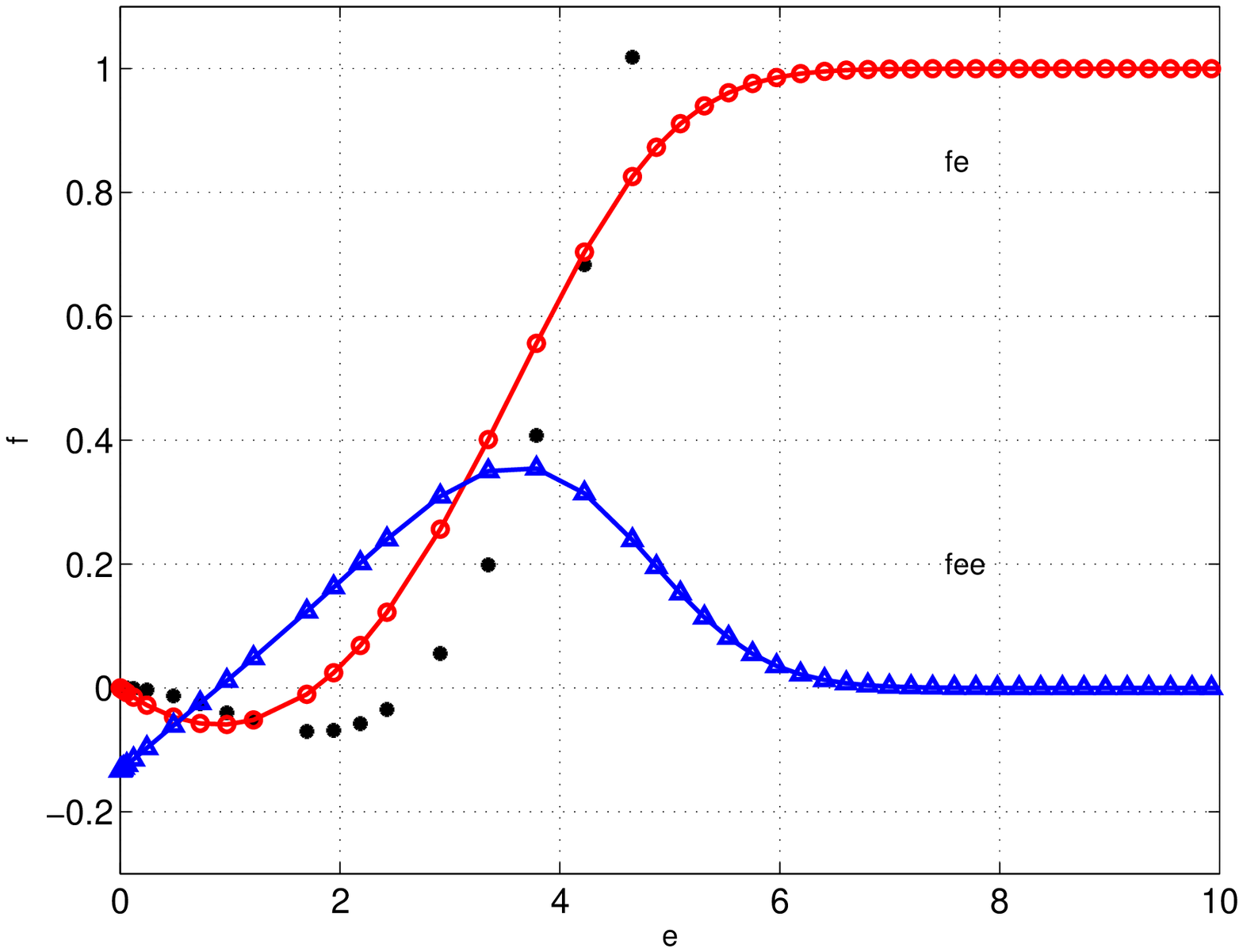}
\caption{Normal and reverse flow solutions to Falkner-Skan model for $\beta =-1.5$.
The symbols $\bullet$ denote values of $f(\eta)$.} 
	\label{fig:Falkner}
\end{figure}

Figure~\ref{fig:Falkner} shows the results of the two numerical solutions for a different value 
of $\beta$, namely $\beta = -0.15$.
In the top frame, we have the normal flow and in the bottom frame, we display the reverse flow solution.
In both cases the solutions were computed by introducing a truncated boundary and solving the IVP in the starred variables on $\eta^* \in [0, 20]$ with $h_0^*=1$, $h_1^* = 5$ in the top frame and
$h_0^*=15$, $h_1^* = 25$ in the bottom frame.
In this case, we achieved convergence of the numerical results within eight and seven iterations, respectively.
For the sake of clarity, we omit to plot the solutions in the starred variables computed during the iterations.
Moreover, we display the solutions only $\eta \in [0, 10]$.

As far as the reverse flow solutions are concerned, in table~\ref{tab:Compare} we compare the missing initial condition computed by the ITM for several values of $\beta$ with results available in the literature.
\begin{table}[!htb]
\caption{Comparison for the reverse flow skin-friction coefficients ${\displaystyle \frac{d^2f}{d\eta^{2}}(0)}$.
For all cases we used $h_0^* = 15$ and $h_1^* = 25$.
The iterations were, from top to bottom line: $8$, $7$, $9$, $7$, and $7$.}
\vspace{.5cm}
\renewcommand\arraystretch{1.3}
	\centering
		\begin{tabular}{r@{.}lcccc}
\hline 
\multicolumn{2}{c}%
{$\beta$} &
Stewartson \cite{Stewartson:1954:FSF} &
Asaithambi \cite{Asaithambi:1997:NMS} & Auteri et al. \cite{Auteri:2012:GLS}&
ITM \\[1ex]
\hline
$-0$ & $025$ & $-0.074$ & & $$ & $-0.074366$ \\
$-0$ & $05$  & $-0.108$ & & $$ & $-0.108271$ \\
$-0$ & $1$    & $-0.141$ & $-0.140546$ & $-0.140546$ & $-0.140546$ \\
$-0$ & $15$  & $-0.132$ & $-0.133421$ & $-0.133421$ & $-0.133421$ \\
$-0$ & $18$  & $-0.097$ & $-0.097692$ & $-0.097692$ & $-0.097692$ \\
\hline			
		\end{tabular}
	\label{tab:Compare}
\end{table}
The agreement is really very good.
It is remarkable that among the studies quoted in the introduction only a few report data related to the reverse flow solutions. 

\begin{figure}[!htb]
	\centering
\psfrag{b}[][]{$ \beta $} 
\psfrag{Blasius}[l][l]{Blasius} 
\psfrag{fee}[][]{${\displaystyle \frac{d^2f}{d\eta^2}} $} 
\includegraphics[width=.7\textwidth]{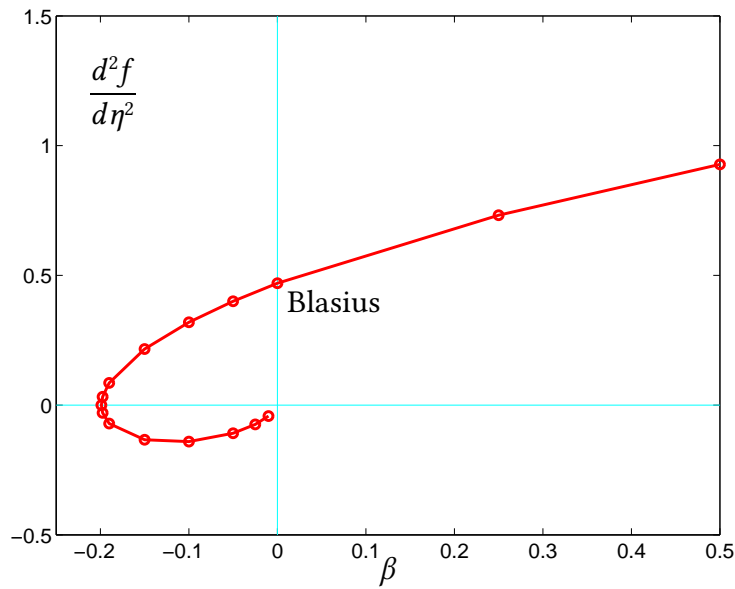}
\caption{Missing initial conditions to Falkner-Skan model for several values of $\beta$.
Positive values determine normal flow, and instead, negative values define reverse flow solutions.} 
	\label{fig:bfee}
\end{figure}
In figure~\ref{fig:bfee} we plot the behaviour of missed initial condition versus $\beta$.
The solution found by the data in table~\ref{tab:Itera2} is plotted in this figure, but not the one found in table~\ref{tab:Itera1} because this is very close to the Blasius solution.
A good initial choice of the initial iterates of $h^*$, for a given value of $\beta$, is obtained by employing values close to the one used in a successful attempt made for a close value of $\beta$.
It is interesting to note that, for values of $\beta < \beta_{\min}$ the ITM continued to iterate endlessly, whatever a couple of starting values for $h^*$ are selected.

Our extended algorithm has shown a kind of robustness because it is able to get convergence even when, for a chosen value of $h^*$, the IVP solver stops before arriving at the selected truncated boundary getting a wrong value of $\Gamma(h^*) = -1$.
On the other hand, the secant method gives an overflow error when this happens for two successive iterates of $h^*$.

\begin{figure}[!htb]
	\centering
\psfrag{e}[][]{\small $ \eta $} 
\psfrag{f}[][]{\small $f(\eta)$} 
\psfrag{fe}[][]{$ {\displaystyle \frac{df}{d\eta}} $} 
\psfrag{fee}[][]{$ {\displaystyle \frac{d^2f}{d\eta^2}} $} 
\includegraphics[width=.7\textwidth]{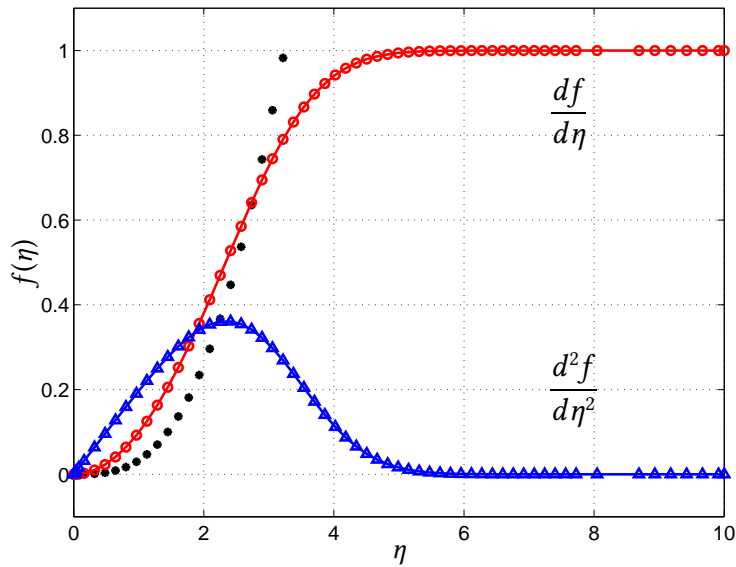}
\caption{Numerical solutions to Falkner-Skan model for $\beta = -0.1988376$. We notice that $ {\displaystyle \frac{d^2f}{d\eta^2}}(0) = 0$ and values of $f(\eta)$ are marked by $\bullet$.} 
	\label{fig:FSHL}
\end{figure}

The value of $\beta_{\min}$, corresponding to a separation point at $\eta=0$, can be found by the ITM by considering $\beta$ as a continuation parameter.
In figure~\ref{fig:FSHL} we plot the unique solution for the limiting value $\beta_{\min}$, where $\beta_{\min}$ is given by equation (\ref{eq:bmin}).
Let us notice that this is a normal flow solution. 

The results reported so far have been found by a variable order adaptive multi-step IVP solver that was coupled with the simple secant method.
The adaptive solver uses a relative and absolute error tolerance, for each component of the numerical solution, both equal to $10\mbox{D}-6$. 
As well known, the secant method is convergent provided that we select two initial iterates sufficiently close to the root, and its convergence is super-linear with an order of convergence equal to $(1+\sqrt{5})/2$.
As far as a termination criterion for the secant method is concerned, we enforced the conditions
\begin{equation}\label{eq:TC}
| \Gamma(h_j^*) | \le \mbox{Tol} \qquad \mbox{and} \qquad |h_j^*-h_{j-1}^*| \le \mbox{TolR}|h_j^*|+ \mbox{TolA} \ ,
\end{equation}
with $ \mbox{Tol}=\mbox{TolR}=\mbox{TolA}=1\mbox{D}-06$.
 
For other problems, in boundary layer theory, admitting more that one solution or none, depending on the value of a parameter involved, see \cite{Fazio:2009:NTM} or \cite{Fazio:2013:BPF}.

\section{Final remarks and conclusions}
In this paper, we define an ITM 
to provide numerical evidence for the existence and uniqueness of a solution of BVPs defined on infinite intervals on the basis of the theorem proved in section 3. 
The leading point is to establish whether the introduced transformation function has only one real zero or not.
From a numerical viewpoint, several computations may be necessary to understand and to characterize the behaviour of $ \Gamma (\cdot) $.
For the proposed numerical test the transformation function was calculated at test-points while at some points a root-finding method was used. 
In the process we tried to bracket the roots of $ \Gamma (\cdot) = 0 $.
In this way, it was possible to obtain numerical results that are in very good agreement with the values available in the literature.
By setting different values of 
$ \delta$ and $ \sigma $ we get a different transformation function and consequently we must bracket its zeros once again.

Of course, the asymptotic boundary condition is not easily used numerically and for the sake of simplicity, we have chosen to replace it with the same condition at a suitable truncated boundary.
A recent successful way to deal with such an issue is to reformulate the considered problem as a free BVP \cite{Fazio:1992:BPF,Fazio:1994:FSE,Fazio:1996:NAN}; for a survey on this topic see \cite{Fazio:2002:SFB}. 
Recently, Zhang and Chen \cite{Zhang:2009:IMS} have used a free boundary formulation to compute the normal flow solutions of the Falkner-Skan model in the full range $ \beta_{\min} < \beta \le 40$.
They applied a modified Newton's method to compute both the initial velocity and the free boundary.  
A recent approach, proposed by Fazio and Jannelli \cite{Fazio:2014:QUG,Fazio:2014:FDS} see also \cite{Fazio:2017:BII,Fazio:2018:NSG}, used to enforce the asymptotic boundary condition exactly, is to apply a quasi-uniform grid, that has the last mesh point at infinity, and suitable non-standard finite difference schemes

In conclusion, the main results of the present investigation are the following: we have defined a theoretical basis for the ITM;
within the definition of the method we emphasized a numerical test for the existence and uniqueness of the solutions of BVPs defined on infinite intervals;
the ITM computes numerical results that are in very good agreement with the correct values.

\vskip 1truecm

\noindent
{\bf Acknowledgement.} This work was supported by the University of Messina and partially by
the Italian GNCS of INDAM.

\end{document}